\documentclass[11pt]{article}
\usepackage{stmaryrd}
\usepackage{mathrsfs}
\usepackage[centertags]{amsmath}
\usepackage{amsfonts, dsfont}
\usepackage{amssymb}
\usepackage{amsthm}

\theoremstyle{plain}
\newtheorem{thm}{Theorem}[section]
\newtheorem{cor}[thm]{Corollary}
\newtheorem{lem}[thm]{Lemma}
\newtheorem{prop}[thm]{Proposition}

\theoremstyle{definition}
\newtheorem{defn}[thm]{Definition}
\newtheorem{remark}[thm]{Remark}
\newtheorem*{ack}{Acknowledgments}

\newcommand{\bd}{\begin{defn}}
\newcommand{\ed}{\end{defn}}
\newcommand{\bl}{\begin{lem}}
\newcommand{\el}{\end{lem}}
\newcommand{\bp}{\begin{prop}}
\newcommand{\ep}{\end{prop}}
\newcommand{\bt}{\begin{thm}}
\newcommand{\et}{\end{thm}}
\newcommand{\bc}{\begin{cor}}
\newcommand{\ec}{\end{cor}}
\newcommand{\br}{\begin{remark}}
\newcommand{\er}{\end{remark}}
\newcommand{\bdi}{\begin{diagram}}
\newcommand{\edi}{\end{diagram}}
\newcommand{\beq}{\begin{eqn}}
\newcommand{\eeq}{\end{eqn}}
\newcommand{\ba}{\begin{array}}
\newcommand{\ea}{\end{array}}
\newcommand{\bpf}{\begin{proof}}
\newcommand{\epf}{\end{proof}}

\newcommand{\Z}{\mathds{Z}}
\newcommand{\Q}{\mathds{Q}}
\newcommand{\Zp}{\mathds{Z}_{p}}

\newcommand{\be}{\beta}

\newcommand{\Ga}{\Gamma}

\newcommand{\Kc}{\mathcal{K}}
\newcommand{\Hi}{H_{\mathrm{Iw}}}

\newcommand{\m}{\mathfrak{m}}

\DeclareMathOperator{\rank}{rank} \DeclareMathOperator{\Gal}{Gal}
\DeclareMathOperator{\Hom}{Hom} \DeclareMathOperator{\Ext}{Ext}

\newcommand{\plim}{\displaystyle \mathop{\varprojlim}\limits}

\newcommand{\cyc}{\mathrm{cyc}}

\newcommand{\lra}{\longrightarrow}

\newcommand{\sbs}{\subseteq}

\newcommand{\ps}[1]{\llbracket #1 \rrbracket}

\begin{document}

\title{On the homology of Iwasawa cohomology groups}
\author{Meng Fai Lim\footnote{Department of Mathematics, University of Toronto, 40 St. George St.,
Toronto, Ontario, Canada M5S 2E4}}
\date{}
\maketitle

\begin{abstract} \footnotesize
\noindent A fundamental observation of Iwasawa gives a criterion for
a $\Zp\ps{\Ga}$-module to be torsion. In this paper, we study and
prove a certain extension of this criterion which will be applied to
investigate the structure of the homology of the ``Iwasawa
cohomology groups'' of a Galois module over a $p$-adic Lie
extension. Namely, we will answer partially a question of Hachimori
on the structure of the homology of the first Iwasawa cohomology
groups as modules over Iwasawa algebras and study the pseudo-nullity
of the homology of the second Iwasawa cohomology groups and dual
fine Selmer groups.


\end{abstract}

\section{Introduction}

Let $p$ be a fixed prime number. Denote $\Ga$ to be the compact
abelian multiplicative group isomorphic to the additive group $\Zp$.
A well-known criterion of Iwasawa states that a finitely generated
$\Zp\ps{\Ga}$-module $M$ with the property that $M_{\Ga}$ is finite
is necessarily a torsion $\Zp\ps{\Ga}$-module. This criterion is an
easy consequence of the structure theory of finitely generated
$\Zp\ps{\Ga}$-module (see \cite{Ho} for another approach in deriving
the criterion) and has been a useful tool in the study of modules
over the classical Iwasawa algebra. This criterion has been extended
to solvable uniform pro-$p$ group in \cite{BHo, HaO}.  Building on
the ideas of \cite{BHo, HaO}, the author established a relative
version of this generalized criterion and applied it to the study of
pseudo-nullity of fine Selmer groups in \cite{Lim}. In this paper,
we show a sharper version of this result under a certain stronger
assumption (see Proposition \ref{algebra main}). The proof of this
result makes additional use of a technical proposition of Kato
\cite{Ka2}.

Following the classical situation, we will apply this criterion to
study the structure of certain modules which now bears a module
structure over a larger (possibly non-commutative) Iwasawa algebra.
Namely, we will study an extension of the questions of Hachimori
\cite{Ha} on the structure of the homology of the first Iwasawa
cohomology groups which were also studied in \cite{Ka2}. Hachimori's
motivation of studying these questions lies in developing a
criterion to determine the existence of the Euler characteristic of
the fine Selmer group of an elliptic curve and he was able to answer
his question in the case when the Iwasawa cohomology groups are
considered over a $p$-adic Lie extension of dimension 2. In this
paper, we will formulate a refined version of his question and
provide an answer to it for more general Galois modules and $p$-adic
extensions (see Theorem \ref{main}). We are also able to establish
the existence of the Euler characteristic of the fine Selmer group
of an abelian variety for this class of $p$-adic Lie extensions (see
Corollary \ref{finite euler char Z}).

We will also formulate and study certain variant of the question of
Hachimori over an admissible $p$-adic extension. As we will see, the
answer to this question in this form is a consequence of various
conjectures on the structure of the second Iwasawa cohomology groups
(or dual fine Selmer groups) which have been studied in \cite{CS,
JhS, Lim}. We finally note that our criterion (Proposition
\ref{algebra main}) has a pseudo-null analog which in fact follows
as a corollary. This latter criterion will be applied to establish
the pseudo-nullity of certain homology groups of the Iwasawa
cohomology groups (and the dual fine Selmer groups).

We now give a brief description of the layout of the paper. In
Section \ref{ranks of Iwasawa modules}, we will discuss and prove
our algebraic criterion. In Section \ref{Hachimori question}, we
formulate a refinement of Hachimori's question and answer the
question for certain classes of $p$-adic Lie extensions. We will
also establish the existence of the Euler characteristic of the dual
fine Selmer group of an abelian variety over the same class of
$p$-adic Lie extensions.  We then formulate and study a variant of
Hachimori's question over an admissible $p$-adic Lie extension in
Section \ref{Hachimori question2}. In Section 5, we will establish
the pseudo-nullity of certain homology groups of the second Iwasawa
cohomology groups and the dual fine Selmer groups.

\section{Ranks of Iwasawa modules} \label{ranks of Iwasawa modules}

Denote $p$ to be a fixed prime. Let $R$ be a commutative complete
regular local ring with a finite residue field of characteristic
$p$. For a compact pro-$p$ $p$-adic Lie group $G$ without
$p$-torsion, we denote $R\ps{G}$ to be the completed group algebra
which is defined by $\plim R[G/U]$, where the inverse limit is taken
over all open normal subgroups of $G$. It is well known that
$R\llbracket G\rrbracket$ is an Auslander regular ring (cf.
\cite[Theorems 3.26 and 3.30(ii)]{V1} or \cite[Theorem A1(c)]{Lim}).
In particular, the ring $R\ps{G}$ is Noetherian local and has finite
projective dimension. Therefore, it follows that every finitely
generated $R\ps{G}$-module admits a finite free resolution of finite
length. In the case that either the ring $R$ has characteristic zero
or $G$ is a uniform pro-$p$ group, the ring $R\ps{G}$ has no zero
divisors (cf. \cite[Theorem A1(b)]{Lim}), and therefore admits a
skew field $Q(G)$ which is flat over $R\ps{G}$ (see \cite[Chapters 6
and 10]{GW} or \cite[Chapter 4, \S 9 and \S 10]{Lam}). In this case,
if $M$ is a finitely generated $R\ps{G}$-module, we define the
$R\ps{G}$-rank of $M$ to be
$$ \rank_{R\ps{G}}M  = \dim_{Q(G)} Q(G)\otimes_{R\ps{G}}M.$$
For a general compact pro-$p$ $p$-adic Lie group $G$ without
$p$-torsion, the $R\ps{G}$-rank of a finitely generated
$R\ps{G}$-module $M$ is defined by
 $$ \rank_{R\ps{G}} M =
 \displaystyle\frac{\rank_{R\ps{G_0}}M}{|G:G_0|}, $$
where $G_0$ is an open normal uniform pro-$p$ subgroup of $G$. This
is integral and independent of the choice of $G_0$ (cf.
\cite[Section 4]{Lim}). In fact, one has the following formula of
Howson (cf. \cite[Theorem 1.1]{Ho} or \cite[Lemmas 4.3 and
4.5]{Lim}).

\bl \label{relative rank} Let $G$ be a compact pro-$p$ $p$-adic Lie
group without $p$-torsion. Let $M$ be a finitely generated
$R\ps{G}$-module. Then we have the equality
 $$  \rank_{R\ps{G}}M =  \displaystyle\sum_{i\geq 0}
(-1)^i\rank_{R}H_i(G,M).
     $$
 Furthermore, if $N$ is a closed normal subgroup of
$G$ such that $G/N$ has no $p$-torsion, we have
 $$ \rank_{R\ps{G}}M = \displaystyle\sum_{i\geq 0} (-1)^i\rank_{R\ps{G/N}}H_i(N,M)
     $$
 \el

 We say that
$M$ is a \textit{torsion} $R\ps{G}$-module if $\rank_{R\ps{G}} M =
0$. The following result was stated in \cite[Theorem 4.6]{Lim},
although the essences of the proof were essentially in \cite{BHo,
HaO, Ho}.

\bp
 Let $G$ be a compact pro-$p$ $p$-adic Lie
group without $p$-torsion. Suppose that $N$ is a closed normal
subgroup of $G$ such that $N$ is solvable uniform and $G/N$ has no
$p$-torsion. Let $M$ be a finitely generated $R\ps{G}$-module. Then
$M_N$ is a torsion $R\ps{G/N}$-module if and only if $M$ is a
torsion $R\ps{G}$-module and
 $\displaystyle\sum_{n\geq 1}
(-1)^n\rank_{R\ps{G/N}}H_n(N,M) =0$. \ep

The above cited proposition is a natural refinement of the results
in \cite{BHo, HaO} and, of course, the classical case when $R=\Zp$
and $N=G=\Zp$. The proposition is false if one removes the
``solvable'' assumption (see \cite{BHo, Ho} for discussions and
counterexamples on this issue). By imposing a stronger condition on
$H$ and $G$, we can prove that each of the terms in the alternating
sum is actually zero.

\bp \label{algebra main}  Let $G$ be a compact pro-$p$ $p$-adic Lie
group without $p$-torsion. Let $N$ be a closed normal subgroup of
$G$ such that $G/N$ has no $p$-torsion. Suppose that there is a
finite family of closed normal subgroups $N_i$ $(0\leq i\leq r)$ of
$G$ such that $1=N_0\sbs N_1 \sbs \cdots\sbs N_r =N$,
$N_i/N_{i-1}\cong \Zp$ for $1\leq i\leq r$ and such that the action
of $G$ on $N_i/N_{i-1}$ by inner automorphism is given by a
homomorphism $\chi_i:G/N\lra \Zp^{\times}$. Let $M$ be a finitely
generated $R\ps{G}$-module. Then the following statements are
equivalent.
 \vspace{-.05in}
\begin{enumerate}
 \item[$(a)$] $M_N$ is a torsion $R\ps{G/N}$-module.

\vspace{-.1in}
\item[$(b)$]  $M$ is a torsion $R\ps{G}$-module, and $H_n(N,M)$ is a torsion
$R\ps{G/N}$-module for every $n\geq 1$.
\end{enumerate}
 \ep

To prove Proposition \ref{algebra main}, we need a proposition of
Kato \cite{Ka2}. Before stating this proposition, we introduce some
terminology. Let $\mathcal{G}$ be a compact pro-$p$, $p$-adic Lie
group. For a continuous character $\kappa : \mathcal{G} \lra
R^{\times}$ and an $R\ps{\mathcal{G}}$-module $M$, we let
$M(\kappa)$ denote the $R$-module $M$ with the new commuting
$\mathcal{G}$-action given by the twist of the original by $\kappa$.
Since $R$ is a $\Zp$-algebra, there is a natural ring homomorphism
$\Zp\lra R$ which in turn induces a group homomorphism
$\Zp^{\times}\lra R^{\times}$. For any continuous character $\chi :
\mathcal{G} \lra \Zp^{\times}$, by abuse of notation, we denote
$\chi : \mathcal{G} \lra R^{\times}$ to be continuous character
obtained by the composition of the above map $\Zp^{\times}\lra
R^{\times}$ with $\chi$. Then $M(\chi)$ is the
$R\ps{\mathcal{G}}$-module defined as above. We can now state the
following proposition of Kato \cite{Ka2}.

\bp \label{Kato} Let $G$ be a compact $p$-adic Lie group and let $N$
be a closed normal subgroup of $G$. Assume that we are given a
finite family of closed normal subgroups $N_i$ $(0\leq i\leq r)$ of
$G$ such that $1=N_0\sbs N_1 \sbs \cdots\sbs N_r =N$,
$N_i/N_{i-1}\cong \Zp$ for $1\leq i\leq r$ and such that the action
of $G$ on $N_i/N_{i-1}$ by inner automorphism is given by a
homomorphism $\chi_i:G/N\lra \Zp^{\times}$.

Let $M$ be a finitely generated $R\ps{G}$-module, and let $M'$ be a
$R\ps{G}$-subquotient of $M$. Let $m\geq 0$. Then there exists a
finite family $(S_i)_{1\leq i\leq k} = (S_{i,m})_{1\leq i\leq k}$ of
$R\ps{G/N}$-submodules of $H_m(N,M')$ satisfying the following
properties. \vspace{-0.05in}
\begin{enumerate}
\item[$(i)$] $0 = S_0\sbs S_1\sbs \cdots\sbs S_k = H_m(N,M')$.

\vspace{-0.05in}
\item[$(ii)$] For each $j$ $(1\leq j\leq k)$, there is a
$R\ps{G/N}$-subquotient $K= K_j$ of $H_0(N,M)$ and a family
$(s_j(i))_{1\leq i\leq r}= (s(i))_{1\leq i\leq r}$ of nonnegative
integers such that $|\{i |s(i)>0\}|\geq m$ and such that
$S_j/S_{j-1}$ is isomorphic to the twist $K(\prod_{1\leq i\leq
r}\chi_i^{s(i)})$ of $K$. \end{enumerate} \ep

\bpf The proof of \cite[Proposition 4.2]{Ka2} basically carries over
word-for-word.  One only needs to make the cosmetic change by
replacing $\Zp$ in the original proof by $R$. \epf

We will require another lemma. For a continuous character $\kappa :
\mathcal{G} \lra R^{\times}$, we denote $\kappa^{-1}$ to be the
character $\mathcal{G} \lra R^{\times}$ given by $g\mapsto
\kappa(g)^{-1}$.

\bl \label{twist lemma}
 Let $\mathcal{G}$ be a compact pro-$p$ $p$-adic Lie
group without $p$-torsion. Suppose that we are given a continuous
character $\kappa : \mathcal{G} \lra R^{\times}$ and a finitely
generated $R\ps{\mathcal{G}}$-module $M$. If $M$ is a torsion
$R\ps{\mathcal{G}}$-module, so is $M(\kappa)$. \el

\bpf By \cite[Lemma 4.4]{Lim}, $M$ is a torsion
$R\ps{\mathcal{G}}$-module if and only if
$\Hom_{R\ps{\mathcal{G}}}(M, R\ps{\mathcal{G}})=0$. It is a
straightforward exercise to see that
 $$\Hom_{R\ps{\mathcal{G}}}(M(\kappa), R\ps{\mathcal{G}})=
\Hom_{R\ps{\mathcal{G}}}(M, R\ps{\mathcal{G}})(\kappa^{-1}).$$  The
conclusion of the lemma is now immediate. \epf

We can now prove Proposition \ref{algebra main}.

\bpf[Proof of Proposition \ref{algebra main}] It follows from the
above proposition of Kato that $H_m(N,M)$ is a successive extension
of twists of $R\ps{G/H}$-subquotients of $H_0(N,M)$. Since $H_0(N,M)
= M_N$ is a torsion $R\ps{G/N}$-module by hypothesis, so is every
$R\ps{G/N}$-subquotient of $M_N$. By Lemma \ref{twist lemma}, every
twist of a $R\ps{G/N}$-subquotient of $M_N$ is also a torsion
$R\ps{G/N}$-module. The conclusion of Proposition \ref{algebra main}
is now immediate. \epf

\medskip
We end the section giving a variant of Proposition \ref{algebra
main} for pseudo-null modules. This will be used in Section
\ref{pseudo-null section}. Recall that a finitely generated torsion
$R\ps{G}$-module is said to be \textit{pseudo-null} if
$\Ext_{R\ps{G}}^{1}(M, R\ps{G}) = 0$.

\bp \label{pseudo-null main} Let $N$ and $H$ be closed normal
subgroups of $G$ such that $N\sbs H$, $G/H\cong \Zp$ and $G/N$ is a
pro-$p$ group of dimension $\geq 2$ without $p$-torsion. Suppose
further that there is a finite family of closed normal subgroups
$N_i$ $(0\leq i\leq r)$ of $G$
 such that $1=N_0\sbs N_1 \sbs \cdots\sbs N_r =N$, $N_i/N_{i-1}\cong
\Zp$ for $1\leq i\leq r$ and such that the action of $H$ on
$N_i/N_{i-1}$ by inner automorphism is given by a homomorphism
$\chi_i: H/N\lra \Zp^{\times}$. Let $M$ be a finitely generated
$R\ps{G}$-module which is also finitely generated over $R\ps{H}$.
Then the following statements are equivalent.

\vspace{-0.05in}
 \begin{enumerate}
 \item[$(a)$] $M_N$ is a pseudo-null $R\ps{G/N}$-module.

\vspace{-0.05in}
 \item[$(b)$] $M$ is a pseudo-null $R\ps{G}$-module, and
$H_n(N,M)$ is a pseudo-null $R\ps{G/N}$-module for every $n\geq 1$.
 \end{enumerate}
\ep

\bpf
 Since $M$ is a
$R\ps{G}$-module which is finitely generated over $R\ps{H}$, it
follows from a result of Venjakob (cf. \cite[Example 2.3 and
Proposition 5.4]{V2} or \cite[Lemma 5.1]{Lim}) that $M$ is a
pseudo-null $R\ps{G}$-module if and only if $M$ is a torsion
$R\ps{H}$-module. The conclusion of the proposition is now an
immediate consequence of Proposition \ref{algebra main}. \epf

\section{On the questions of Hachimori} \label{Hachimori question}

In this section, we will apply the results obtained in the preceding
section to study the questions of Hachimori. We begin by introducing
some terminology and notation that we shall use throughout this
section.

As before, $p$ will denote a fixed prime. Let $F$ be a number field.
If $p=2$, we assume further that $F$ has no real primes. Let $S$ be
a finite set of primes of $F$ that contains the primes above $p$ and
the infinite primes. We then denote $F_S$ to be the maximal
algebraic extension of $F$ unramified outside $S$. For any algebraic
(possibly infinite) extension $\mathcal{L}$ of $F$ contained in
$F_S$, we write $G_S(\mathcal{L}) = \Gal(F_S/\mathcal{L})$. Let
$F_{\infty}$ be a $p$-adic extension of $F$ which is unramified
outside $S$. We shall assume further that $G= \Gal(F_{\infty}/F)$ is
a compact pro-$p$ Lie group with no $p$-torsion.

Let $R$ be a complete regular local ring with a finite residue field
of characteristic $p$. Let $T$ denote a finitely generated free
$R$-module with a continuous $R$-linear $G_S(F)$-action. (Here the
topology on $T$ is given by its filtration by powers of the maximal
ideal of $R$.) We define the $i$th ``Iwasawa cohomology groups'' (of
$T$ over $F_{\infty}$) to be
 \[ H_{\mathrm{Iw},S}^i(F_{\infty}/F, T) = \plim_L H^i(G_S(L), T), \]
where the inverse limit is taken over all the finite extensions $L$
of $F$ contained in $F_{\infty}$ and with respect to the
corestriction maps. For ease of notation, we will drop the `$S$'. It
is well-known that $\Hi^i(F_{\infty}/F, T)$ is a finitely generated
(compact) $R\ps{G}$-module for every $i$. In fact, it is not
difficult to see that $\Hi^i(F_{\infty}/F, T)=0$ for $i\neq 1,2$. We
now record the following lemma which will be crucial in our
discussion.

\bl \label{spectral sequence} Assume that $H^0(G_S(F), T)=0$. We
then have isomorphisms
  \[ \ba{c} \Hi^2(F_{\infty}/F,T)_G \cong H^2(G_S(F),T)
  ~\mbox{and} \vspace{0.05in}\\
  H_n(G, \Hi^1(F_{\infty}/F,T)) \cong H_{n+2}(G,
 \Hi^2(F_{\infty}/F,T)) ~\mbox{for $n \geq 1$}.  \ea\] \el

\bpf
 These are immediate consequences from analyzing the homological spectral
sequence
$$ H_{i}(G, \Hi^{-j}(F_{\infty}/F,T)) \Longrightarrow H^{-i-j}(G_S(F), T) $$
which can be found, for instance, in either \cite[Proposition
1.6.5(iii)]{FK} or \cite[Theorem 3.1.8]{LS}. \epf

In this section, we are interested in the following question which
was raised and studied in \cite{Ha}.

\medskip \noindent \textbf{Question:} Do the following statements hold?

 (A) $\rank_{R\ps{G}} \Hi^1(F_{\infty}/F, T) = \rank_{R} \Hi^1(F_{\infty}/F, T)_G$.

 (B) $H_n\big(G, \Hi^1(F_{\infty}/F, T)\big)$ is a torsion $R$-module for each $n\geq 1$.

\medskip
It follows from an application of the formula of Howson (cf. Lemma
\ref{relative rank}) that (B) implies (A). Throughout the article,
we will focus on the validity of (B). In his paper \cite{Ha},
Hachimori studied the question in the case when $T$ is the Tate
module of an elliptic curve. In \cite{Ka2}, the question is studied
in the case when $T=\Zp(1)$. We can now state the main result of
this section which gives an affirmative answer to the above question
for certain classes of $p$-adic Lie extensions.

\bt \label{main} Let $F_{\infty}$ be a $p$-adic Lie extension of $F$
with Galois group $G$. Assume that $H^0(G_S(F), T) =0$. Suppose
further that either one of the following statements holds.

\vspace{-0.05in}
 \begin{enumerate}
 \item[$(i)$]  The group $G$ has dimension $\leq 2$ and has no $p$-torsion.

\vspace{-0.05in}
 \item[$(ii)$] There is a finite family of closed normal
subgroups $G_i$ $(0\leq i\leq r)$ of $G$
 such that $1=G_0\sbs G_1 \sbs \cdots\sbs G_r =G$ and $G_i/G_{i-1}\cong
 \Zp$ for every $i$, and $H^2(G_S(F), T)$ is a torsion $R$-module.
\end{enumerate}
\vspace{-0.05in} Then statement $(B)$ holds. \et

Before proving the theorem, we make a remark.

\br The assumption $H^0(G_S(F), T) =0$ is known to hold in many
interesting cases (for instance, $\Zp(r)$ where $r\neq 0$, and the
Tate module of an abelian variety). When $T$ is the Tate module of
an elliptic curve $E$, the theorem under statement (i) was proved in
\cite[Theorem 1.3]{Ha} under further assumptions that
$E(F_{\infty})[p^{\infty}]$ is finite and $H^2(G_S(F_{\infty}),
E[p^{\infty}])=0$. Therefore, our theorem in this case can be viewed
as a refinement of that, namely both in terms of the Galois module
$T$ and that we can do without the two assumptions above. We also
note that if $G$ has dimension $\leq 2$, then $G$ will satisfy the
group theoretical assumption of statement (ii). (This is clear if
$G$ is abelian or has dimension 1. In the case when $G$ is
nonabelian of dimension 2, this follows from \cite[Proposition
7.1]{GSK}.) The point of (i) is that if the dimension of $G$ is less
than or equal to 2, one does not require the torsionness condition
on $H^2(G_S(F), T)$ in (ii).  \er

\bpf[Proof of Theorem \ref{main}]

(i) By the assumption that $\dim G\leq 2$ and the second isomorphism
in Lemma \ref{spectral sequence}, we have $H_n(G,
\Hi^1(F_{\infty}/F,T)) = 0$ for all $n\geq 1$. In particular,
statement (B) holds.

(ii) It follows from Lemma \ref{spectral sequence} that $H_0\big(G,
\Hi^2(F_{\infty}/F, T)\big)$ is isomorphic to $H^2(G_S(F),T)$ which
is $R$-torsion by the hypothesis. By an application of Proposition
\ref{algebra main}, we then have that $H_n\big(G,
\Hi^2(F_{\infty}/F, T)\big)$ is $R$-torsion for all $n\geq 1$. By
the second isomorphism in Lemma \ref{spectral sequence}, this in
turn implies that $H_n\big(G, \Hi^1(F_{\infty}/F, T)\big)$ is
$R$-torsion for all $n\geq 1$.  Thus, statement (B) holds. \epf

We now define the dual fine Selmer group. Let $v$ be a prime in $S$.
For every finite extension $L$ of $F$ contained in $F_S$, we define
 \[K_v^2(T/L) = \bigoplus_{w|v}H^2(L_w, T),\]
where $w$ runs over the (finite) set of primes of $L$ above $v$. If
$\mathcal{L}$ is an infinite extension of $F$ contained in $F_S$, we
define
\[\Kc_v^2(T/\mathcal{L}) = \plim_L K_v^2(T/L),\]
where the inverse limit is taken over all finite extensions $L$ of
$F$ contained in $\mathcal{L}$. For any algebraic extension
$\mathcal{L}$ of $F$ contained in $F_S$, the \textit{dual fine
Selmer group} of $T$ over $\mathcal{L}$ (with respect to $S$) is
defined to be
\[ Y(T/\mathcal{L}) = \ker\Big(\Hi^2(\mathcal{L}/F, T)\lra
\bigoplus_{v\in S}\Kc_v^2(T/\mathcal{L})
\Big). \]

Since $Y(T/F_{\infty})$ is a $R\ps{G}$-submodule of
$\Hi^2(F_{\infty}/F,T)$, the following is an immediate consequence
of Proposition \ref{algebra main} and Lemma \ref{spectral sequence}.

\bt \label{finite euler char} Let $F_{\infty}$ be a $p$-adic Lie
extension of $F$ with Galois group $G$. Assume that $H^0(G_S(F), T)
=0$. Suppose that there is a finite family of closed normal
subgroups $G_i$ $(0\leq i\leq r)$ of $G$
 such that $1=G_0\sbs G_1 \sbs \cdots\sbs G_r =G$ and $G_i/G_{i-1}\cong
 \Zp$ for $1\leq i\leq r$, and that $H^2(G_S(F), T)$ is $R$-torsion.
Then $H_n(G, Y(T/F_{\infty}))$ is $R$-torsion for every $n\geq 0$.
\et

For the remainder of the section, we will assume that our ring $R$
is $\Zp$. As mentioned in the introduction, one of the motivation of
studying the question of Hachimori lies in developing a criterion
for the existence of the Euler characteristic of the dual fine
Selmer group of a $p$-adic representation (see \cite[Theorem
3.1]{Ha}). Following \cite{Ha}, we say that the $G$-Euler
characteristic of $Y(T/F_{\infty})$ is defined if $H_n(G,
Y(T/F_{\infty}))$ is finite for all $n$. (Note that since we are
working with the dual fine Selmer group, we have to consider
$G$-homology instead of $G$-cohomology. However, the $G$-homology of
the dual fine Selmer group of $T$ is dual to the $G$-cohomology of
the fine Selmer group of $T^{\vee}(1)$. Therefore, our definition
here is compatible with that in \cite{Ha}.) An interesting
by-product of our discussion is that we can establish the existence
of the Euler characteristic of the dual fine Selmer group.

\bc \label{finite euler char Z} Let $A$ be an abelian variety over
$F$. Denote $T$ to be the Tate module of the dual abelian variety
$A^*$ of $A$. Let $F_{\infty}$ be a $p$-adic Lie extension of $F$
with Galois group $G$. Suppose that there is a finite family of
closed normal subgroups $G_i$ $(0\leq i\leq r)$ of $G$
 such that $1=G_0\sbs G_1 \sbs \cdots\sbs G_r =G$ and $G_i/G_{i-1}\cong
 \Zp$ for every $i$. Then if $Y(T/F)$ is finite,
the $G$-Euler characteristic of $Y(T/F_{\infty})$ is defined. \ec

\bpf
 Since $A^*(F)[p^{\infty}]$ is finite, we have $H^0(G_S(F), T) =0$.
 On the other hand, for each $v\in S$, it follows easily from Mattuck's theorem \cite{Mat} that
 $A(F_v)[p^{\infty}]$ is finite. By Tate local duality, this in turn implies that $H^2(F_v, T)$ is finite. Therefore,
 $Y(T/F)$ is finite if and only if $H^2(G_S(F), T)$ is finite. Thus, the
 hypotheses of Theorem \ref{finite euler char} are satisfied and the required
 conclusion follows.  \epf

\br To obtain a similar conclusion of Corollary \ref{finite euler
char Z} for a general module $T$, one needs to replace the condition
of $Y(T/F)$ being finite by the condition that $H^2(G_S(F),T)$ is
finite. This is because for a general $T$, one may not have the
finiteness of $H^2(F_v, T)$. \er

\section{On a variant of the question of Hachimori} \label{Hachimori question2}

In this section, we will study a variant of the question of
Hachimori. In particular, we show that this question is a
consequence of various conjectures on the structure of the second
Iwasawa cohomology groups.

As before, $p$ will denote a fixed prime. Let $F$ be a number field,
where we assume further that it has no real primes when $p=2$.  Let
$S$ be a finite set of primes of $F$ which contains the primes above
$p$ and the infinite primes. Following \cite{CS}, we say that
$F_{\infty}$ is an $S$-admissible $p$-adic Lie extension of $F$ if
(i) $\Gal(F_{\infty}/F)$ is a compact pro-$p$ $p$-adic Lie group
without $p$-torsion, (ii) $F_{\infty}$ contains the cyclotomic
$\Zp$-extension $F^{\cyc}$ of $F$ and (iii) $F_{\infty}$ is
contained in $F_S$. Write $G=\Gal(F_{\infty}/F)$,
$H=\Gal(F_{\infty}/F^{\cyc})$ and $\Ga=\Gal(F^{\cyc}/F)$. As before,
we denote $T$ to be a finitely generated free $R$-module with a
continuous $R$-linear $G_S(F)$-action, where $R$ is a complete
regular local ring with a finite residue field of characteristic
$p$. We will like to study the following variant of the question of
Hachimori.

\medskip \noindent \textbf{Question:} Do the following statement hold?

 $(\beta)$ $H_n\big(H, \Hi^1(F_{\infty}/F, T)\big)$ is a torsion
 $R\ps{\Ga}$-module for each $n\geq 1$.

\medskip
To facilitate further discussion, we record the following lemma (cf.
\cite[Lemma 3.2]{CS} or \cite[Lemma 5.2]{Lim}).

\bl \label{fg La H} Let $F_{\infty}$ be an $S$-admissible $p$-adic
Lie extension of $F$. Then the following statements are equivalent.

\smallskip
$(a)$ $Y(T/F^{\mathrm{cyc}})$ is a finitely generated $R$-module.

\smallskip
$(b)$ $\Hi^2(F^{\cyc}/F, T)$ is a finitely generated $R$-module.

\smallskip
$(c)$ $Y(T/F_{\infty})$ is a finitely generated $R\ps{H}$-module.

\smallskip
$(d)$  $\Hi^2(F_{\infty}/F, T)$ is a finitely generated
$R\ps{H}$-module.\el

We now recall the following conjecture which has been studied in
\cite{CS, JhS, Lim}.

\medskip \noindent \textbf{Conjecture A:} For
any number field $F$,\, one (and hence all) of the statements in
Lemma \ref{fg La H} holds.

\bp \label{finitely generated} Let $F_{\infty}$ be an $S$-admissible
$p$-adic Lie extension of $F$ with Galois group $G$. Assume that
Conjecture A holds. Then $H_n(H, \Hi^j(F_{\infty}/F, T))$ is a
finitely generated $R$-module for every $n\geq 1$ and $j=1,2$.

In particular, statement $(\be)$ holds. \ep

Before proving the proposition, we record the following variant of
Lemma \ref{spectral sequence} which will be required in the proof of
the proposition and subsequent part of this paper.

\bl \label{spectral sequence2} Let $F_{\infty}'$ be an infinite
$p$-adic extension of $F$ contained in $F_{\infty}$ which has the
property that $\Gal(F_{\infty}'/F)$ has no $p$-torsion. Write
$N=\Gal(F_{\infty}/F_{\infty}')$. We then have isomorphisms
  \[\ba{c} \Hi^2(F_{\infty}/F,T)_N \cong \Hi^2(F_{\infty}'/F,T)~\mbox{and} \vspace{0.05in}\\
   H_n(N, \Hi^1(F_{\infty}/F,T)) \cong H_{n+2}(N,
 \Hi^2(F_{\infty}/F,T)) ~\mbox{for $n \geq 1$}. \ea \]
\el

\bpf
 Since $F_{\infty}'$ is
an infinite pro-$p$ extension of $F$, we have
$\Hi^0(F_{\infty}'/F,T)=0$. The proof then proceeds as in Lemma
\ref{spectral sequence}. \epf

\bpf[Proof of Proposition \ref{finitely generated}] Since we are
assuming that Conjecture A is valid, we have that
$\Hi^2(F_{\infty}/F,T)$ is a finitely generated $R\ps{H}$-module. It
is then an easy exercise (or see \cite[Lemma 3.2.3]{LS}) to show
that $H_{n}(H, \Hi^2(F_{\infty}/F,T))$ is a finitely generated
$R$-module for every $n\geq 1$. By the second isomorphism in Lemma
\ref{spectral sequence2}, this in turn implies that $H_{n}(H,
\Hi^1(F_{\infty}/F,T))$ is a finitely generated $R$-module for every
$n\geq 1$. In particular, statement $(\be)$ holds. \epf

When $T=\Zp(1)$, Conjecture A is precisely the classical conjecture
of Iwasawa which asserts that the $\mu$-invariant of
$\Gal(K(F^{\cyc})/F^{\cyc})$ vanishes (see \cite{Iw}). Here
$K(F^{\cyc})$ is the maximal unramified pro-$p$ extension of
$F^{\cyc}$ where every prime of $\mathcal{L}$ above $p$ splits
completely. We shall call this conjecture the Iwasawa
$\mu$-conjecture for $F^{\cyc}$. It is interesting to note that the
general statement of Conjecture A turns out to be a consequence of
this classical conjecture of Iwasawa (cf. \cite[Theorem 3.1]{Lim}).
However, even in this classical setting, the conjecture is only
proved in the case when $F$ is abelian over $\Q$ (see \cite{FW,
Sin}). Building on this, we at least have the following
unconditional validity of statement $(\be)$.

\bc \label{finitely generated 2} Let $F_{\infty}$ be an
$S$-admissible $p$-adic Lie extension of $F$ with Galois group $G$.
Assume further that $F(\mu_{2p}, T/\m T)$ is a finite $p$-extension
of an abelian extension $F'$ of $\Q$. Then statement $(\be)$ holds.
\ec

\bpf
 It follows from an application of the main theorem of \cite{FW, Sin} and \cite[Theorem 3]{Iw2} that
the Iwasawa $\mu$-conjecture holds for $F(\mu_{2p}, T/\m T)^{\cyc}$.
By \cite[Theorem 3.5]{Lim}, this implies that Conjecture $A$ holds
for $T$ over $F(\mu_{2p}, T/\m T)^{\cyc}$. By \cite[Lemma 3.2]{Lim},
this in turn implies that Conjecture $A$ holds for $T$ over
$F^{\cyc}$. The conclusion of the corollary now follows from
Proposition \ref{finitely generated}. \epf

In particular, Conjecture $A$ holds for elliptic curves over $\Q$
that have a rational $p$-isogeny (cf. \cite[Corollary 3,6]{CS}).
Therefore, we can apply Corollary \ref{finitely generated 2} to
obtain $(\be)$ unconditionally. We shall give some examples that are
not an elliptic curve. For the rest of this paragraph, we shall
assume that $p \geq 5$. Let $n\geq 1$. Let $a, b ,c$ be integers
such that $1\leq a, b,c < p^n$, $a+b+c = p^n$ and at least one of
(and hence at least two of) $a, b, c$ is not divisible by $p$. Let
$J=J_{a,b,c}$ be the Jacobian variety of the curve $y^{p^n} =
x^a(1-x)^b$. Set $T$ to be the Tate module of the $p$-division
points of $J$. Write $R = \Zp[\zeta_{p^n}]$ and $\pi =
1-\zeta_{p^n}$. Note that $\pi$ is a generator of the maximal ideal
of $R$. It follows from \cite[Chap.\ II \S\S4, Chap.\ II Theorem
5B]{Ih} (see also \cite[Section 2]{Gr} for the case $n=1$) that $T$
is a free $R$-module of rank one and that $G_S(\Q(\mu_{p^n}))$ acts
trivially on $T/\pi T$. Therefore, we may apply the preceding
corollary to obtain validity of $(\be)$ for any admissible $p$-adic
Lie extension of $\Q(\mu_p)$.

\medskip
We will now investigate the statement $(\be)$ under the following
weaker assumption on $\Hi^2(F_{\infty}/F, T)$.

\medskip \noindent \textbf{Conjecture WL$_{/F_{\infty}}$:} \textit{For
any $S$-admissible extension $F_{\infty}$ of $F$,\,
$\Hi^2(F_{\infty}/F, T)$ is a torsion $R\ps{G}$-module.}

\medskip We note that Conjecture WL will follow from Conjecture A.
However unlike Conjecture A, Conjecture WL$_{/F_{\infty}}$ does not
have good descent properties and hence the dependence on
$F_{\infty}$. In general, one can show that if Conjecture
WL$_{/F_{\infty}}$ holds, then Conjecture WL$_{/L_{\infty}}$ holds
for any solvable extension $L_{\infty}$ of $F_{\infty}$ (cf.\
\cite[Proposition 7.2]{Lim}). Despite this, Conjecture
WL$_{/F_{\infty}}$ is known to hold in many cases. For instance, it
is well-known that Conjecture WL$_{/F_{\infty}}$ holds for every
$F_{\infty}$ if $T=\Zp(1)$ (cf. \cite[Theorem 5]{Iw} and
\cite[Theorem 6.1]{OcV}). In the case when $E$ is an elliptic curve
defined over $\Q$ with good reduction at $p$ and $F$ is an abelian
extension of $\Q$, a deep theorem of Kato \cite[Theorem 12.4]{Ka1}
asserts that Conjecture WL$_{/F^{\cyc}}$ is valid.

As it turns out, we are not able to establish an analogue statement
as Proposition \ref{finitely generated} assuming Conjecture WL only.
However, we do have the following partial result in this direction
whose proof, which we omit, is similar to Theorem \ref{main}.

\bt \label{torsion descent} Let $F_{\infty}$ be an $S$-admissible
$p$-adic Lie extension of $F$ with Galois group $G$. As before,
denote $H= \Gal(F_{\infty}/F^{\cyc})$. Suppose that one of the the
following statements holds. \vspace{-0.05in}
 \begin{enumerate}
 \item[$(i)$]  The group $H$ has dimension $\leq 2$ and has no $p$-torsion.

\vspace{-0.05in}
 \item[$(ii)$] Suppose that there is a finite family
of closed normal subgroups $H_i$ $(0\leq i\leq r)$ of $G$
 such that $1=H_0\sbs H_1 \sbs \cdots\sbs H_r =H$, $H_i/H_{i-1}\cong
\Zp$ for $1\leq i\leq r$ and such that the action of $G$ on
$H_i/H_{i-1}$ by inner automorphism is given by a homomorphism
$\chi_i:G/H\lra \Zp^{\times}$. Also, suppose that Conjecture
WL$_{/F^{\cyc}}$ holds.
\end{enumerate}
Then statement $(\be)$ holds. \et

Combining the preceding theorem with Kato's theorem \cite[Theorem
12.4]{Ka1} on the validity of Conjecture WL$_{/F^{\cyc}}$, we have
the following corollary.

\bc \label{torsion descent2} Set $F= \Q(\mu_p, a_1,..., a_n)$, where
each $a_i$ is a nonzero element of the maximal abelian extension of
$\Q$ which is not a root of unity. Let $F_{\infty} =
\Q(\mu_{p^{\infty}}, a_1^{-p^{\infty}},...,a_n^{-p^{\infty}})$.
Suppose that $T$ is one of the following Galois modules.
\vspace{-0.05in}
 \begin{enumerate}
 \item[$(i)$]  The Tate module of an elliptic curve over $\Q$
with good reduction at $p$. Here $T$ is a module over $\Zp$.

\vspace{-0.05in}
 \item[$(ii)$] A Galois invariant lattice of the Galois
 representation associated to a normalized cuspidal eigenform $f$ of weight $\geq 2$,
 tame level $N$ and character $\psi$ that is ordinary at
 $p$. Here $T$ is a module over the ring of integers of the completion of the number
field $K_f$ at some prime above $p$, where $K_f$ is the number field
generated by the Fourier coefficients of $f$ and the values of
$\psi$.

\vspace{-0.05in}
 \item[$(iii)$] The Galois representation associated to a nearly ordinary Hida
 deformation $($cf. $\cite[P.181 (4)]{JhS})$ satisfying the conditions
 $(\textbf{Nor})$ and $(\textbf{Irr})$ as given in $\cite[P.
 182]{JhS}$. Here $T$ is a module over $\mathcal{O}\ps{X}$.

\end{enumerate} Then statement $(\be)$ holds. \ec

\bpf
 As is well-known (and noting that $F$ is an abelian extension of $\Q$), the validity of Conjecture
WL$_{/F^{\cyc}}$ in cases (i) and (ii) follows from \cite[Theorem
12.4]{Ka1}. The validity of Conjecture WL$_{/F^{\cyc}}$ in case
(iii) follows from that in (i) and (ii) by a standard control
argument (or see \cite[Proposition 7.3]{Lim}). The conclusion of the
corollary then follows from an application of Theorem \ref{torsion
descent}. \epf

\br We like to mention that in a recent paper of Aribam \cite{A}, he
was able to establish Conjecture A for certain elliptic curves and
modular forms (which are defined over $\Q$) over the cyclotomic
$\Zp$-extension of certain number fields $F$ that are \textit{non
abelian} over $\Q$. His examples can therefore be applied to give
(unconditional) examples of validity of statement $(\be)$ that are
not covered by our Corollary \ref{finitely generated 2} and
Corollary \ref{torsion descent2}. \er

\section{Pseudo-nullity of homology of Iwasawa cohomology groups} \label{pseudo-null section}

We retain the notation of Section \ref{Hachimori question2}. The aim
of this section is to prove the following theorem which refines
\cite[Theorem 5.7]{Lim} under a stronger assumption.

\bt \label{pseudo-null descent} Let $F_{\infty}$ be an
$S$-admissible $p$-adic Lie extension of $F$. Assume that Conjecture
A holds. Suppose that $F_{\infty}'$ is another $S$-admissible
$p$-adic Lie extension of $F$ which satisfies the following
properties.

\vspace{-0.05in}
 \begin{enumerate}
 \item[$(i)$] $F'_{\infty}$ is contained in $F_{\infty}$, and we
 write $N=\Gal(F_{\infty}/F_{\infty}')$.

\vspace{-0.05in}
 \item[$(ii)$] For each $v\in S$, the decomposition
group of $\Gal(F_{\infty}'/F)$ at $v$ has dimension $\geq 2$.

\vspace{-0.05in}
 \item[$(iii)$] There is a finite family of closed normal subgroups $N_i$
$(0\leq i\leq r)$ of $\Gal(F_{\infty}/F)$ such that $1=N_0\sbs N_1
\sbs \cdots\sbs N_r =N$, $N_i/N_{i-1}\cong \Zp$ for $1\leq i\leq r$
and such that the action of $H$ on $N_i/N_{i-1}$ by inner
automorphism is given by a homomorphism $\chi_i:G/N\lra
\Zp^{\times}$.

\end{enumerate}
Then $Y(T/F_{\infty}')$ is a pseudo-null
$R\ps{\Gal(F_{\infty}'/F)}$-module if and only if $Y(T/F_{\infty})$
is a pseudo-null $R\ps{\Gal(F_{\infty}/F)}$-module and $H_n\big(N,
\Hi^2(F_{\infty}/F, T)\big)$ is a pseudo-null
$R\ps{\Gal(F_{\infty}'/F)}$-module for every $n\geq 1$. \et

\bpf
 By assumption (ii) and \cite[Lemma 5.3]{Lim}, $Y(T/F_{\infty}')$ is a pseudo-null
$R\ps{\Gal(F_{\infty}'/F)}$-module if and only if
$\Hi^2(F_{\infty}'/F, T)$ is a pseudo-null
$R\ps{\Gal(F_{\infty}'/F)}$-module. We have a similar statement for
$Y(T/F_{\infty})$ and $\Hi^2(F_{\infty}/F, T)$. The conclusion of
the theorem now follows from an application of Proposition
\ref{pseudo-null main} and Lemma \ref{spectral sequence2}. \epf

Combining the above theorem with Lemma \ref{spectral sequence2}, we
have the following corollary which will answer a pseudo-null analog
of statement $(\be)$.

\bc \label{pseudo-null descent corollary} Retain the assumptions of
Theorem \ref{pseudo-null descent}. Suppose that $Y(T/F_{\infty}')$
is a pseudo-null $R\ps{\Gal(F_{\infty}'/F)}$-module. Then
$H_n\big(N, \Hi^1(F_{\infty}/F, T)\big)$ is a pseudo-null
$R\ps{\Gal(F_{\infty}'/F)}$-module for every $n\geq 1$. \ec

We end by giving an example. Let $E$ be the elliptic curve $79A1$ of
Cremona's table which is given by
\[ y^2 + xy + y = x^3 + x^2 -2x.\]
Take $p=3$ and $F = \Q(\mu_3)$. Let $S$ be the set of primes of $F$
lying above $3, 79$ and $\infty$. Let $T$ denote either the Tate
module of $E$ or the Galois representation attached to the Hida
family associated to the weight 2 newform corresponding to $E$. It
was shown in \cite{Jh} that $Y(T/F_{\infty}')$ is a pseudo-null
$R\ps{\Gal(F_{\infty}'/F)}$-module when $F_{\infty}' =
\Q(\mu_{3^{\infty}}, 79^{-3^{\infty}})$. Since $\Q(\mu_{3^{\infty}},
79^{-3^{\infty}})$ satisfies statement (b) of Theorem
\ref{pseudo-null descent} (cf.\ \cite[Lemma 3.9]{HV}), one can apply
Theorem \ref{pseudo-null descent} to conclude that $Y(T/F_{\infty})$
is a pseudo-null $R\ps{\Gal(F_{\infty}/F)}$-module,
$H_n\big(\Gal(F_{\infty}/F'_{\infty}), \Hi^i(F_{\infty}/F,
T_5E)\big)$ is a pseudo-null $R\ps{\Gal(F'_{\infty}/F)}$-module for
every $n\geq 1$ and $i=1,2$, and
$H_n\big(\Gal(F_{\infty}/F'_{\infty}), Y(T_5E/F_{\infty})\big)$ is a
pseudo-null $R\ps{\Gal(F'_{\infty}/F)}$-module for every $n\geq 1$,
where $F_{\infty}$ is one of the following $S$-admissible $3$-adic
Lie extensions:
\[ \Q(\mu_{3^{\infty}},
3^{3^{-\infty}}, 79^{3^{-\infty}}), \quad
L_{\infty}(79^{3^{-\infty}}),\quad L_{\infty}(3^{3^{-\infty}},
79^{3^{-\infty}}). \] Here $L_{\infty}$ is the unique
$\Z_3^2$-extension of $F$.

\begin{ack}
     This work was written up when the author is a Postdoctoral fellow at the GANITA Lab
    at the University of Toronto. He would like to acknowledge the
    hospitality and conducive working conditions provided by the GANITA
    Lab and the University of Toronto.
        \end{ack}

\end{document}